\documentclass{article}
\usepackage{latexsym,amssymb}
\usepackage{graphicx}
\usepackage[noadjust]{cite}
\input amssym.def \input amssym

\newtheorem{te}{Theorem}[section]

\newtheorem{de}[te]{Definition}
\newtheorem{lm}[te]{Lemma}

\newtheorem{pp}[te]{Proposition}

\newtheorem{ex}[te]{Example}
\newtheorem{qu}[te]{Question}
\def\dokaz{\noindent{\bf Proof. }}
\def\kraj{\hfill $\Box$ \par \vspace*{2mm} }

\def\widemid{\hspace{1mm}\widetilde{\mid}\hspace{1mm}}

\newcommand{\zve}[1]{{{}^*\hspace{-0.5mm}#1}}
\newcommand{\zvez}[1]{{{}^*\hspace{-1mm}#1}}

\newcommand{\tac}[1]{{{}^\bullet\hspace{-0.5mm}#1}}
\def\zs{\forall}
\def\po{\exists}
\def\str{\rightarrow}
\def\Str{\Rightarrow}
\def\rtS{\Leftarrow}
\def\dl{\Leftrightarrow}
\def\gstr{\hspace{-0.1cm}\uparrow}
\def\dstr{\hspace{-0.1cm}\downarrow}
\def\ps{\subseteq}

\def\rest{\upharpoonright}

\def\cU{{\cal U}}
\def\cV{{\cal V}}
\def\cW{{\cal W}}

\def\cF{{\cal F}}
\def\cG{{\cal G}}
\def\cH{{\cal H}}
\def\cM{{\cal M}}
\def\cP{{\cal P}}
\def\cQ{{\cal Q}}
\def\cR{{\cal R}}
\def\cC{{\cal C}}
\def\bN{{\mathbb{N}}}
\def\bZ{{\mathbb{Z}}}

\begin{document}
\begin{center}
           {\huge \bf Congruence of ultrafilters}
\end{center}
\begin{center}
{\small \bf Boris  \v Sobot}\\[2mm]
{\small  Department of Mathematics and Informatics, University of Novi Sad,\\
Trg Dositeja Obradovi\'ca 4, 21000 Novi Sad, Serbia\\
e-mail: sobot@dmi.uns.ac.rs\\
ORCID: 0000-0002-4848-0678}
\end{center}
\begin{abstract}
We continue the research of the relation $\widemid$ on the set $\beta \bN$ of ultrafilters on $\bN$, defined as an extension of the divisibility relation. It is a quasiorder, so we see it as an order on the set of $=_\sim$-equivalence classes, where $\cF=_\sim\cG$ means that $\cF$ and $\cG$ are mutually $\widemid$-divisible. Here we introduce a new tool: a relation of congruence modulo an ultrafilter. We first recall the congruence of ultrafilters modulo an integer and show that $=_\sim$-equivalent ultrafilters do not necessarily have the same residue modulo $m\in \bN$. Then we generalize this relation to congruence modulo an ultrafilter in a natural way. After that, using iterated nonstandard extensions, we introduce a stronger relation, which has nicer properties with respect to addition and multiplication of ultrafilters. Finally, we also introduce a strengthening of $\widemid$ and show that it also behaves well in relation to the congruence relation.
\end{abstract}


\section{Introduction}

Let $\bN$ be the set of natural numbers. The relation $\widemid$, an extension of the divisibility relation $\mid$ on $\bN$ to the set $\beta \bN$ of ultrafilters on $\bN$, was introduced in \cite{So1} and further investigated in \cite{So2, So3, So4, So5}. The main idea was to understand the impact of various properties of $\mid$ to $\widemid$ and possibly, learning about the $\widemid$-hierarchy, to acquire better understanding of $\mid$. In this paper we will make another step in that direction, considering possible extensions of the congruence relations to $\beta \bN$ and their relation to $\widemid$, as well as to the operations of addition and multiplication on $\beta \bN$.\\

When working with the set of ultrafilters $\beta S$ on a set $S$ it is common to identify each element $s\in S$ with the principal ultrafilter $\{A\subseteq S:s\in A\}$. Having that in mind, any binary operation $\star$ on $S$ can be extended to $\beta S$ as follows: for $A\subseteq S$,
\begin{equation}\label{eqn1}
A\in\cF\star\cG\dl\{s\in S:s^{-1}A\in\cG\}\in\cF,
\end{equation}
where $s^{-1}A=\{t\in S:s\star t\in A\}$. If $(S,\star)$ is a semigroup equiped with the discrete topology, $(\beta S,\star)$ becomes a compact Hausdorff right-topological semigroup. The base sets for the topology are (clopen) sets $\overline{A}=\{\cF\in\beta S:A\in\cF\}$. Many aspects of structures obtained in this way were examined in \cite{HS}.\\

Every function $f:\bN\str \bN$ can be extended uniquely to a continuous $\widetilde{f}:\beta \bN\str\beta \bN$: the ultrafilter $\widetilde{f}(\cF)$ is generated by $\{f[A]:A\in\cF\}$. This was used in \cite{So1} to define analogously an extension of a binary relation $\rho$ on $\bN$ to a relation $\widetilde{\rho}$ on $\beta \bN$: $\cF\widetilde{\rho}\cG$ if and only if for every $A\in\cF$ the set $\rho[A]:=\{n\in \bN:(\po a\in A)a\rho n\}$ is in $\cG$. This coincides with the so-called canonical way of extending relations from $\bN$ to $\beta \bN$ described in \cite{Gor}. It turned out that the extension $\widemid$ of the divisibility relation $\mid$ has a simple equivalent definition, more convenient for practical use:
$$\cF\widemid\cG\dl\cF\cap\cU\subseteq\cG,$$
where $\cU=\{A\in P(\bN)\setminus\{\emptyset\}:A\gstr=A\}$ is the family of sets upward closed for $\mid$. $\widemid$ is a quasiorder, so we think of it as an order on the set of $=_\sim$-equivalence classes, where $\cF=_\sim\cG$ if and only if $\cF\widemid\cG$ and $\cG\widemid\cF$. We say that $C\subseteq \bN$ is {\it convex} if for all $x,y\in C$ and all $z$ such that $x\mid z$ and $z\mid y$ holds $z\in C$. All ultrafilters from the same $=_\sim$-equivalence class $\cC$ have the same convex sets. Clearly, each equivalence class $\cC$ is determined by $\cF\cap\cU$ (for any $\cF\in\cC$), or by the family of convex sets belonging to any $\cF\in\cC$. 

An ultrafilter $\cF$ is divisible by some $n\in \bN$ if and only if $n\bN:=\{nk:k\in\bN\}\in\cF$. If $\cF\in\bN$ as well, $n\widemid\cF$ holds if and only if $n\mid\cF$. Hence, we can write just $n\mid\cF$ in case $n\in\bN$.

Especially useful are prime ultrafilters $\cP$: those $\widemid$-divisible only by 1 and themselves. This is equivalent to $P\in\cP$, where $P$ is the set of prime numbers.

The $\widemid$-hierarchy can be naturally divided into two parts. The ``lower" part, $L$, can be divided into levels: $L=\bigcup_{l<\omega}\overline{L_l}$, where $L_l=\{p_1p_2\dots p_l:p_1,p_2,\dots,p_l\mbox{ are prime}\}$ is the set of natural numbers having exactly $l$ (not necessarily distinct) prime factors. Some nice properties of $L$ were established in \cite{So3}; for example every ultrafilter in $\overline{L_l}$ has exactly $l$ prime ingredients (but being divisible by the $n$-th power of a prime $\cP$ is not the same as being divisible by $\cP$ $n$ times). The ``upper" part, however, is much more complicated. It contains the maximal $=_\sim$-class, $MAX$, consisting of ultrafilters divisible by all $n\in \bN$, and consequently by all $\cF\in\beta \bN$ (\cite{So4}, Lemma 4.6). Another interesting class is $NMAX$, maximal among $\bN$-{\it free} ultrafiters (those that are not divisible by any $n\in \bN$), see \cite{So5}, Theorem 5.4. A set belonging to an $\bN$-free ultrafilter is called an $\bN$-{\it free} set.\\

The paper is organized as follows. In Section 2 several well-known results of elementary number theory are employed to obtain results about the congruence of ultrafilters modulo an integer in connection with the divisibility relation $\widemid$. In Section 3 we recapitulate basic definitions about $\omega$-hyperextensions, obtained by iterating nonstandard extensions of the set $\bZ$. Tensor pairs play an important role here. They were first considered by Puritz in \cite{P2}; Di Nasso proved several useful characterizations and coined the term (see \cite{DN}). Most of the results in Section 3 are taken from Luperi Baglini's thesis \cite{L1}, where the concept of a tensor pair is implemented in the surrounding of $\omega$-hyperextensions. In Section 4 we define congruence modulo an ultrafilter and find several conditions equivalent to this definition. The next section deals with a stronger relation, and we prove some results connecting it to $\widemid$ and operations of addition and multiplication of ultrafilters. In Section 6 we define another version of divisibility, obtained in a natural way from the strong congruence relation, and get some basic results about it. The last section contains several remarks and open questions.\\

{\bf Notation.} $\bN$ is the set of natural numbers (without zero), $\omega=\bN\cup\{0\}$, $P$ is the set of prime numbers and $\bZ$ the set of integers. The calligraphic letters $\cF,\cG,\cH,\dots$ are reserved for ultrafilters, and small letters $x,y,z,\dots$ for integers (both standard and nonstandard). For $A\subseteq \bN$, $A\gstr=\{n\in \bN:\po a\in A\;a\mid n\}$ and $A\dstr=\{n\in \bN:\po a\in A\;n\mid a\}$. If $m,r\in \bN$, then $\bZ_m=\{0,1,\dots,m-1\}$ and $mA+r=\{mn+r:n\in A\}$. Finally, $\cU=\{A\in P(\bN)\setminus\{\emptyset\}:A\gstr=A\}$ and $\cV=\{A\in P(\bN)\setminus\{\bN\}:A\dstr=A\}$.

Because we use $\zve\bN$ for a nonstandard extension of $\bN$, to avoid confusion we will not denote $\beta\bN\setminus\bN$ with $\bN^*$. Likewise, we will avoid writing $A^2$ for $A\times A$, since this notation had another meaning in papers preceding this one.

\section{Congruence modulo integer}\label{kongm}

Let $m\in \bN$ and let $\bZ_m$ be given the discrete topology. The homomorphism $h_m:\bN\str \bZ_m$ is defined as follows: $h_m(n)$ is the residue of $n$ modulo $m$. $h_m$ extends uniquely to a continuous function $\widetilde{h_m}:\beta \bN\str \bZ_m$. The next results follows directly from \cite{HS}, Corollary 4.22.

\begin{pp}\label{hom}
$h_m$ is a homomorphism, both for addition and multiplication of ultrafilters.
\end{pp}

As described in the Introduction, the relation $\equiv_m$ of congruence modulo $m$ can be extended to a relation $\widetilde{\equiv_m}$ on $\beta \bN$: $\cF\widetilde{\equiv_m}\cG$ if and only if, for every $A\in\cF$, $\{n\in\bN:(\po a\in A)n\equiv_ma\}\in\cG$. Recall that the kernel of a function $h:\bN\str \bN$ is the relation ${\rm ker}h=\{(x,y)\in \bN\times\bN:h(x)=h(y)\}$.

\begin{pp}\label{kernel}
(\cite{So1}, Theorem 2.13) If $h:\bN\str \bN$ and $\rho=\ker h$, then $\widetilde{\rho}=\ker\widetilde{h}$.
\end{pp}

Thus, for $m\in \bN$ the extension of $\equiv_m$ to $\beta \bN$ coincides with the definition found in \cite{HS}: $\cF\widetilde{\equiv_m}\cG$ if and only if $h_m(\cF)=h_m(\cG)$. In particular, $r<m$ is the residue of $\cF\in\beta \bN$ modulo $m$ ($\cF\widetilde{\equiv_m}r$) if and only if $m\bN+r\in\cF$. For practical reasons, we will denote the extension of $\equiv_m$ to $\beta \bN$ also by $\equiv_m$ from now on.


The congruence of ultrafilters modulo integer is not new, but it was mostly marginally mentioned; for example the following interesting result has only the status of a comment in \cite{HS}.

\begin{pp}
(\cite{HS}, Comment 11.20) For every $\cF\in\beta \bN$ and every $U\in\cF$ there is a neighborhood $\bar A$ of $\cF$ such that $A\ps U$ and for all $\cG\in{\bar A}\setminus A$ and all $m\in \bN$ holds $\cG\equiv_m\cF$.
\end{pp}

We begin with a simple result about the solvability of a system of congruences in $\beta \bN$. A system such that its every finite subsystem has a solution in $\beta \bN$ will be called {\it feasible}.

\begin{lm}\label{feasys}
(a) Let $x\equiv_{m_i}a_i$ (for $i=0,1,\dots,k$, $a_i\in\bZ$ and $m_i\in\bN$) be a finite system of congruences. It has a solution in $\beta \bN\setminus \bN$ if and only if it has a solution in $\bN$.

(b) The system $x\equiv_{m_i}a_i$ (for $i\in\omega$, $a_i\in\bZ$ and $m_i\in\bN$) of congruences has a solution in $\beta \bN$ if and only if it is feasible.
\end{lm}

\dokaz (a) Let $\cF\in\beta \bN\setminus \bN$ be a solution of the given system. Then $A_i:=\{x\in \bN:x\equiv_{m_i}a_i\}\in\cF$ for each $i=0,1,\dots,k$. Hence $A:=\bigcap_{i=0}^kA_i\in\cF$, and any $x\in A$ is a solution of the given system.

On the other hand, if $s\in \bN$ is a solution and $u=lcm(m_0,m_1,\dots,m_k)$ (the least common multiplier of $m_0,m_1,\dots,m_k$), then all the elements of the set $B=\{x\in \bN:x\equiv_us\}$ are also solutions. Thus every $\cF\in\overline{B}\setminus B$ is a solution of the system in $\beta \bN\setminus \bN$.\\

(b) One direction is trivial, so assume the given system to be feasible. Let $A_i=\{x\in \bN:x\equiv_{m_i}a_i\}$. By the assumption, every finite subsystem of the given system has a solution, so the family $\{\overline{A_i}:i<\omega\}$ has the finite intersection property. Since all the sets $\overline{A_i}$ are closed, it follows that $A=\bigcap_{i<\omega}\overline{A_i}$ is nonempty, and any $\cF\in A$ is a solution of the given system.\kraj

Since $=_\sim$-equivalence classes within $L$ are singletons (\cite{So3}, Corollary 5.10), each class in $L$ trivially contains ultrafilters congruent only to one residue modulo $m$. We want to investigate for which systems of congruences there is a $=_\sim$-class such that all its ultrafilters satisfy  it. Clearly, such a system must be feasible. On the other hand, by Lemma \ref{feasys} a feasible system $S$ has a solution $\cG\in\beta\bN$ so we can assume that it is a system of all congruences satisfied by $\cG$ (we will call such a system {\it maximal}). Also, every congruence $x\equiv_{m_i}r_i$ is equivalent to a system of congruences modulo mutually prime factors of $m_i$, so we can assume that all $m_i$ are powers of primes themselves. Let $Q_S=\{p\in P:\cG\equiv_{p^n}0\mbox{ for all }n\in \bN\}$ and $T_S=P\setminus Q_S$. As a special case, if $T_S=\emptyset$, all ultrafilters from the class $MAX$ satisfy $S$.

$A\subset \bN$ is an {\it antichain} if there are no distinct $a,b\in A$ such that $a\mid b$.

\begin{te}\label{jedanost}
For every maximal feasible system $S$ of congruences $x\equiv_{p^n}r_{p,n}$ (for $n\in\omega$, $p\in P$ and $r_{p,n}<p^n$) such that $T_S$ is infinite there is an $=_\sim$-equivalence class $\cC\not\subseteq L$ such that $\cF\equiv_{p^n}r_{p,n}$ for all $\cF\in\cC$.
\end{te}

\dokaz We consider two cases.\\

$1^\circ$ $Q$ is infinite. Let $\{q_i:i\in\omega\}$ and $\{t_i:i\in\omega\}$ be enumerations of $Q_S$ and $T_S$ respectively. For $i\in\omega$ let $s_i=\min\{n\in\bN:\cG\not\equiv_{t_i^n}0\}$. We construct, by recursion on $n$, a set $A=\{a_n:n\in\omega\}$ such that $a_n<a_{n+1}$ and:

(1) $a_n\in t_i^{s_i+n}\bN+r_{t_i,s_i+n}$ for $i<n$;

(2) $t_n^{s_n}\mid a_n$;

(3) $q_j^n\mid a_n$ for every $j<n$.

Start with choosing any $a_0\in t_0^{s_0}\bN$. Assume that $a_n$ is constructed. We want to choose $a_{n+1}$ satisfying the system $x\equiv_{t_i^{s_i+n+1}}r_{t_i,s_i+n+1}$ for $i\leq n$, $x\equiv_{t_{n+1}^{s_{n+1}}}0$ and $x\equiv_{q_j^n}0$ for $j\leq n$. By the Chinese remainder theorem this system has a solution in $\bN$ such that $a_{n+1}>a_n$. Clearly, obtained $a_{n+1}$ satisfies conditions (1)-(3).

$A$ is an antichain: for all $m<n$, $a_m<a_n$ implies that $a_n\nmid a_m$, and $t_m^{s_m}\mid a_m$ and (1) imply that $a_m\nmid a_n$. Let $\cC$ be the $=_\sim$-equivalence class of any ultrafilter containing $A$. Every ultrafilter $\cF\in\cC$ contains $A\gstr$ and $A\dstr$, so it contains $A=A\gstr\cap A\dstr$ as well. Condition (3) clearly implies that $A$ intersects each level $L_l$ only in finitely many elements, so $\cF\notin L$, and in particular $\cF$ is nonprincipal. By (1), $A\setminus(t_i^{s_i+n}\bN+r_{t_i,s_i+n})$ is finite for all $i$ and all $n$, hence $\cF\equiv_{t_i^{s_i+n}}r_{t_i,s_i+n}$. By (3), $\cF\equiv_{q_i^n}0$ for all $i\in\omega$ and $n\in\bN$. Thus $\cF$ satisfies all congruences of the given system.\\

$2^\circ$ $Q$ is finite. We repeat the construction from case $1^\circ$, but for $j\geq |Q_S|$ (when we ``run out" of elements from $Q_S$) instead of $q_j$ in condition (3) we use some elements $t_i\in T$ for $i>n$. (This condition is needed here only to ensure that $\cF\notin L$.)\kraj

\begin{pp}\label{free}
(\cite{So5}, Lemma 5.2) If $A$ is an $\bN$-free set, then $A\not\subseteq n_1\bN\cup n_2\bN\cup\dots\cup n_k\bN$ for any $n_1,n_2,\dots,n_k\in\bN\setminus\{1\}$.
\end{pp}

\begin{ex}
(1) Let us show that the condition of $T_S$ being infinite in the theorem above is necessary. Consider a system $S$ consisting of $x\equiv_{t_i}r_i$ (for some primes $t_0,t_1,\dots,t_{l-1}$ and some nonzero $r_i<t_i$) and $x\equiv_{p^n}0$ for all $p\in P\setminus\{t_0,t_1,\dots,t_{l-1}\}$ and all $n\in\bN$. Let us show that there can be no $=_\sim$-class $\cC$ such that all $\cF\in\cC$ satisfy $S$. Assume the opposite. Then every such $\cF$ contains all sets in $\cU_N:=\{A\in\cU:A\mbox{ is }\bN\mbox{-free}\}$: by Proposition \ref{free} each $A\in\cU_N$ must contain an element $a$ mutually prime to all $t_0,t_1,\dots,t_{l-1}$. Hence $a\mid\cF$ implies $a\bN\in\cF$, and therefore $A\in\cF$. This means that $\cF\cap\cU=\cU_N\cup\{n\bN:n\in\bN\land t_i\nmid n\mbox{ for all }i=0,1,\dots,l-1\}$. But now, if we change any of the $r_i$'s into another nonzero value we stay inside the same class $\cC$.

(2) In the class $NMAX$ of $\widemid$-maximal $\bN$-free ultrafilters one can find an ultrafilter congruent to $r$ modulo $m$ for any $0<r<m$ such that $gcd(m,r)=1$. Namely, the family $\cU_N\cup\{\bN\setminus n\bN:n>1\}\cup\{m\bN+r\}$ has the finite intersection property: for any given $A\in\cU_N$ and $n_0,n_1,\dots,n_k\in\bN\setminus\{1\}$, since $A$ is $\bN$-free, Proposition \ref{free} says that there is $a\in A$ mutually prime to all of $m,n_0,\dots,n_k$. By the Chinese remainder theorem the system $x\equiv_mr$, $x\not\equiv_{n_i}0$, $x\equiv_a0$ has a solution, and it belongs to $A\cap(m\bN+r)\cap\bigcap_{0\leq i\leq k}(\bN\setminus n_i\bN)$.
\end{ex}

Now we will prove a result describing which residues modulo a given prime can appear in the same $=_\sim$-class; first we need the following definition. A set $S$ of residues modulo $p\in P$ is {\it a geometric set of residues} if there are $s$ and $r$ such that $0\leq s<p$, $0<r<p$ and $S=\{rest(sr^k,p):k\in \bN\}$, where $rest(x,p)$ is the residue of $x$ modulo $p$.

\begin{te}\label{geomset}
Let $p\in P$ and let $S\subseteq\{0,1,\dots,p-1\}$. There is an $=_\sim$-equivalence class $\cC$ such that the set of residues of ultrafilters $\cF\in\cC$ is exactly $S$ if and only if $S$ is a geometric set of residues.
\end{te}

\dokaz ($\rtS$) First assume that $S=\{s_0,\dots,s_{l-1}\}$ is a geometric set of residues, where $s_i=rest(s_0r^i,p)$ (for $i=0,1,\dots,l-1$) are exactly all distinct residues of numbers $s_0r^k$ modulo $p$. If $S=\{0\}$, which happens for $s_0=0$, any $=_\sim$-class of ultrafilters divisible by $p$ (i.e.\ containg the set $p\bN$) will do. Otherwise, by Dirichlet's prime number theorem, there are primes $s\equiv_ps_0$ and $b\equiv_pr$. Let $B=\{sb^k:k\in\omega\}$, $\cU'=\{U\in\cU:U\cap B\neq\emptyset\}$ and $\cV'=\{V\in\cV:\bN\setminus V\notin\cU'\}$. Then the family $\cU''=\cU'\cup\cV'$ has the finite intersection property: $\cU'$ is closed for finite intersections, and every $V\in\cV'$ contains $B$. Let $\cC$ be the $=_\sim$-equivalence class determined by $\cU''$. For every $\cF\in\cC$ we have $B\in\cF$ (since $B\cup\{b^k:k\in\omega\}\in\cV'$ and $\bN\setminus\{b^k:k\in\omega\}\in\cU'$) and $B\subseteq\bigcup_{i=0}^{l-1}(p\bN+s_i)$, so every such $\cF$ is congruent to some $s_i$ modulo $p$. On the other hand, for each $i\in\{0,1,\dots,l-1\}$ the family $\cU''\cup\{p\bN+s_i\}$ has the finite intersection property: $B$ contains infinitely many elements from each of the sets $p\bN+s_i$, and finite intersections of sets from $\cU''$ contain all but finitely many elements from $B$, so they also intersect $p\bN+s_i$. Hence there is an ultrafilter $\cF\in\cC$ such that $\cF\equiv_p s_i$.\\

($\Str$) Now assume $S$ is the set of residues modulo $p$ of ultrafilters $\cF\in\cC$ for some $=_\sim$-equivalence class $\cC$. Every singleton is clearly a geometric set of residues (obtained by choosing the quotient $r=1$), so we will assume $|S|>1$. Let $\cW$ be the family of all convex sets belonging to all $\cF\in\cC$. Since the elements of $S$ are all possible residues of ultrafilters $\cF\in\cC$, there is $C\in\cW$ (a finite intersection of sets from $(\cU\cup\cV)\cap\cF$) such that $C\subseteq\bigcup_{k=0}^{l-1}(p\bN+s_k)$ (otherwise $\cW\cup\{\bN\setminus\bigcup_{k=0}^{l-1}(p\bN+s_k)\}$ would have the finite intersection property). 

Let $q$ be a primitive root modulo $p$ (this means that for every $0<r<p$ there is $k\in \bN$ such that $q^k\equiv_pr$; see \cite{Bur} for more details). Let $S=\{s_0,\dots,s_{l-1}\}$, where $s_i=rest(q^{k_i},p)$, $k_0<k_1<\dots<k_{l-1}$ and for each $s_i$ the smallest $k_i$ is chosen. If we denote $r_i=k_i-k_0$ for $0<i<l$, then $s_i=rest(s_0q^{r_i},p)$.\\

\underline{Claim 1.} The set $R:=\{r_i:0<i<l\}$ is closed for the $gcd$ (greatest common divisor) operation.

{\it Proof of Claim 1.} Let $0<i<j<l$. Take $A_0$ to be the set of $\mid$-minimal elements of $C\cap(p\bN+s_0)$. By recursion on $k$, let $A_{3k+1}$ be the set of $\mid$-minimal elements of $C\cap A_{3k}\gstr\cap(p\bN+s_i)$, $A_{3k+2}$ the set of $\mid$-minimal elements of $C\cap A_{3k+1}\gstr\cap(p\bN+s_j)$ and $A_{3k+3}$ the set of $\mid$-minimal elements of $C\cap A_{3k+2}\gstr\cap(p\bN+s_0)$. Each of the sets $A_m$ (for $m\in\omega$) must be nonempty, since otherwise $C\subseteq (C\setminus A_0\gstr)\cup(C\cap A_0\gstr\setminus A_1\gstr)\cup\dots\cup(C\cap A_{m-1}\gstr)$, and each of the (convex) sets on the right would miss one of the sets $p\bN+s_0$, $p\bN+s_i$ or $p\bN+s_j$, so it could not belong to all ultrafilters in $\cC$.

Now let $d=gcd(r_i,r_j)$. By B\' ezout's lemma there are $a',b'\in \bZ$ such that $a'r_i+b'r_j=d$. By replacing $a',b'$ with their residues modulo $p-1$ we get $a,b\in \bZ_{p-1}$ such that $ar_i+br_j\equiv_{p-1}d$. Let $m=3(a+b)$ and let $\langle c_i:0\leq i<m\rangle$ be a $\mid$-chain in $C$ of length $m$ such that $c_i\in A_i$ (it exists since $A_{m-1}\neq\emptyset$). Let $c_{i+1}=c_id_i$; then $d_{3k}\equiv_p q^{r_i}$ and $d_{3k}d_{3k+1}\equiv_p q^{r_j}$ for all $k$. Hence
\begin{eqnarray*}
e &:=& d_0d_3\dots d_{3(a-1)}d_{3a}d_{3a+1}d_{3(a+1)}d_{3(a+1)+1}\dots d_{3(a+b-1)}d_{3(a+b-1)+1}\\
  &\equiv_p& (q^{r_i})^a(q^{r_j})^b=q^{ar_i+br_j}\equiv_p q^d
\end{eqnarray*}
(in the last equality we used Fermat's little theorem). But $c_0e$ is divisible by $c_0$ and divides $c_m$; since $C$ is convex, $c_0e\in C$ and hence $d\in R$.\\

\underline{Claim 2.} $rest(tr_1,p-1)\in R$ for all $t\in \bN$.

{\it Proof of Claim 2} is similar to (though simpler than) the one from Claim 1. We construct a $\mid$-chain $\langle c_i:0\leq i\leq 2t-2\rangle$ such that $c_i\in p\bN+s_0$ for odd $i$ and $c_i\in p\bN+s_1$ for even $i$. If $c_{i+1}=c_id_i$, we get $c_0d_1d_3\dots d_{2t-3}\equiv_p q^{tr_1}$, so $tr_1\equiv_{p-1}r_j$ for some $r_j\in R$.\\

Now, since $r_1<r_2<\dots<r_{l-1}$, the two Claims show that $R$ must have the form $R=\{ir_1:0<i<l\}$. But then $s_i\equiv s_0(q^{r_1})^i$, which is what we wanted to prove.\kraj

\section{$\omega$-hyperextensions of $\bZ$}

In the previous two papers, \cite{So4} and \cite{So5}, we employed nonstandard methods (more precisely, the superstructure approach) to get more information on the relation $\widemid$. We will continue that practice here. However, now we turn to extensions of the set $\bZ$ of all integers instead of $\bN$. The reason is, of course, that we want to use the operation of subtraction. Let $X$ be a set containing a copy of $\bZ$ consisting of atoms: none of the elements of $X$ contains as an element any of the other relevant sets. Let $V_0(X)=X$, $V_{n+1}(X)=V_n(X)\cup P(V_n(X))$ for $n\in\omega$ and $V(X)=\bigcup_{n<\omega}V_n(X)$. $V(X)$ is then called a {\it superstructure}. The rank of an element $x\in V(X)$ is the smallest $n\in\omega$ such that $x\in V_n(X)$.

If $V(X)$ is a superstructure, its {\it nonstandard extension} is a pair $(V(Y),*)$, where $V(Y)$ is a superstructure with the set of atoms $Y$ and $*:V(X)\rightarrow V(Y)$ is a rank-preserving function such that $A\subseteq\zve A$ for $A\subseteq X$, $\bZ\subset\zve \bZ$, $\zve X=Y$ and satisfying the Transfer principle (we delay the formulation of Transfer until later, since we will need a more general version).

A nonstandard extension $(V(Y),*)$ of $V(X)$ is a $\kappa$-{\it enlargement} if for every family $F$ of subsets of some set in $V(X)$ with the finite intersection property such that $|F|<\kappa$ there is an element in $\bigcap_{A\in F}\zve A$. $\kappa$-enlargements are known to exist in ZFC.

For an excellent introduction to nonstandard methods we refer the reader to \cite{G}.

The connection between a nonstandard extension and $\beta \bZ$ is given by the function $v:\zve \bZ\str\beta \bZ$, defined by $v(x)=\{A\subseteq \bZ:x\in\zve A\}$. $v$ is onto whenever $(V(Y),*)$ is a ${\goth c}^+$-enlargement.

\begin{pp}\label{slaganjev}
(\cite{NR}, Lemma 1) For every $x\in\zve Z$ and every $f:\bZ\str\bZ$, $v(\zve f(x))=\widetilde{f}(v(x))$.
\end{pp}

More information about $v$ can be found in \cite{NR} and \cite{L1}. The following proposition is Theorem 3.1 of \cite{So4}, adjusted for extensions of $\bZ$ (instead of $\bN$).

\begin{pp}\label{ekviv}
The following conditions are equivalent for every two ultrafilters $\cF,\cG\in\beta \bZ$:

(i) $\cF\widemid\cG$;

(ii) in every ${\goth c}^+$-enlargement $V(Y)$, there are $x,y\in\zve \bZ$ such that $v(x)=\cF$, $v(y)=\cG$ and $x\zvez\mid y$;

(iii) in some ${\goth c}^+$-enlargement $V(Y)$, there are $x,y\in\zve \bZ$ such that $v(x)=\cF$, $v(y)=\cG$ and $x\zvez\mid y$.


\end{pp}

First, let us establish that we can use all previously obtained results about $\zve \bN$ while working with $\zve \bZ$. In every extension $V(Y)$ the nonstandard set $\zve \bZ$ consists of $\zve \bN$, another (``inverted") copy of $\zve \bN$ (containing negative nonstandard numbers) and zero. For $x,y\in\zve \bZ$, $x\zvez\mid y$ holds if and only if $|x|\;\zvez\mid\;|y|$.

The situation with $\beta \bZ$ is similar. Let, for $A\subseteq \bZ$, $-A:=\{-a:a\in A\}$; likewise, for $\cF\in\beta \bN$ let $-\cF:=\{-A:A\in\cF\}$. Then every ultrafilter in $\beta \bZ$ (except the principal ultrafilter identified with zero) contains either $\bN$ or $-\bN$, so $\beta \bZ=\beta \bN\cup\{-\cF:\cF\in\beta \bN\}\cup\{0\}$. The family $\cU_Z:=\{U\in P(\bZ)\setminus\{\emptyset\}:U\gstr=U\}$ of upward closed subsets of $\bZ$ consists of sets $V\cup -V\cup\{0\}$ for $V\in\cU$, and divisibility in $\beta \bZ$ is naturally defined as: $\cF\widemid\cG$ if and only if $\cF\cap\cU_Z\subseteq\cG$. Thus, $\cF\widemid\cG$ if and only if $|\cF|\widemid|\cG|$ (for absolute values of ultrafilters defined in the obvious way).

We will write $\cF-\cG$ instead of $\cF+(-\cG)$. So $A\in\cF-\cG$ if and only if $\{n\in \bZ:n-A\in\cG\}\in\cF$, where $n-A=\{n-a:a\in A\}$. Note that there can be no confusion with this notation, since $\cF-\cG$ is exactly the ultrafilter obtained by extending the subtraction operation from $\bZ$ to $\beta \bZ$, as defined in (\ref{eqn1}).\\

A nonstandard extension $(V(Y),*)$ of $V(X)$ is called a {\it single superstructure model} if $Y=X$. The existence of such model was proved in \cite{Ben}. In a single superstructure model it is possible to iterate the star-function, since it is defined for all elements in the range of $*$.

\begin{de}
Let $(V(X),*)$ be a single superstructure model with $\bZ\subseteq X$. Define recursively, for $x\in V(X)$, $S_0(x)=x$ and $S_{k+1}(x)=\zve(S_k(x))$ for all $k\in\omega$. For $A\subseteq X$ the set $\tac A=\bigcup_{k<\omega}S_k(A)$ is called an $\omega$-hyperextension of $A$.
\end{de}

Now, any $(V(X),S_k)$ is a nonstandard extension, and $(V(X),\bullet)$ is also a nonstandard extension. Moreover, we have the following.

\begin{pp}\label{prenosnatac}
(\cite{L1}, Proposition 2.5.7) If $(V(X),*)$ is a single superstructure model which is a ${\goth c}^+$-enlargement, then $(V(X),S_k)$ for every $k\in\omega$ and $(V(X),\bullet)$ are also ${\goth c}^+$-enlargements.
\end{pp}

We will call a single superstructure model $(V(X),*)$ which is a ${\goth c}^+$-enlargement a $\omega$-{\it hyperenlargement}.

Now we can use the Transfer principle within any of the mentioned extensions. Recall that a first-order formula $\varphi(x_1,x_2,\dots,x_n)$ is bounded if all its quantifiers are bounded, i.e.\ of the form $(\forall x\in y)$ or $(\exists x\in y)$. In the Transfer principle the free variables $x_1,x_2,\dots,x_n$ that appear in $\varphi(x_1,x_2,\dots,x_n)$ can take values of elements $a_1,a_2,\dots,a_n\in V(X)$ and in $\varphi(\zve a_1,\zve a_2,\dots,\zve a_n)$ they are replaced with their star-counterparts. Any $k$-ary relation $A\in V(X)$ appearing as an atomic subformula in $\varphi$ is also considered like a free variable and gets replaced with $\zve A$.\\

{\it The Transfer principle.} For every bounded formula $\varphi$ and every $a_1,a_2,\dots$, $a_n\in V(X)$, in $V(X)$ $\varphi(a_1,a_2,\dots,a_n)$ holds if and only if $\varphi(S_k(a_1),S_k(a_2),\dots$, $S_k(a_n))$ holds (for any $k\in \bN$) if and only if $\varphi(\tac a_1,\tac a_2,\dots,\tac a_n)$ holds.

As a simple application of Transfer let us show that $\zve(x+y)=\zve x+\zve y$ for $x,y\in\tac Z$, a fact that we will need later. If $z=x+y$, Transfer implies that $\zve z=\zve x+\zve y$. Likewise, $\zve(x\cdot y)=\zve x\cdot\zve y$.

\begin{pp}\label{nivoi}
(\cite{L1}, Proposition 2.5.3)

(a) For $k\leq l$ and $A\subseteq \bZ$, $S_k(A)=S_l(A)\cap S_k(\bZ)$. Consequently, $S_k(A)=\tac A\cap S_k(\bZ)$.

(b) For $h:\bZ\str \bZ$ and $x\in S_k(\bZ)$, $\tac h(x)=S_k(h)(x)$.
\end{pp}

Let us comment on the iterated version of the divisibility relation. It is common to omit $*$ (or, more generally, $S_k$) in formulas in front of the relations $=$ and $\in$ and arithmetical operations, in order to simplify notation. Let us show that it is justified to do the same with the divisibility relation, even when working in an $\omega$-hyperextension. Firstly, $(x,y)\in S_k(\mid)$ can hold only if $x,y\in S_k(\bZ)$. On the other hand, for $x\in S_k(\bN)$, $y\in S_k(\bZ)$ and $l>k$, we will show that $(x,y)\in S_k(\mid)$ if and only if $(x,y)\in S_l(\mid)$.

$(x,y)\in S_k(\mid)$ means that there is $z\in S_k(\bZ)$ such that $y=xz$. But $S_k(\bZ)\subseteq S_l(\bZ)$, so $(x,y)\in S_l(\mid)$ follows. In the other direction, if $(x,y)\in S_l(\mid)$ for some $l>k$, and $y=xz$, then $z\in S_k(\bZ)$ so $(x,y)\in S_k(\mid)$ as well. Thus, there will be no ambiguity if we drop the stars and write simply $x\mid y$ instead of $(x,y)\in S_k(\mid)$.

\begin{de}
For $\cF\in\beta\bZ$, $\mu_n(\cF)=\{x\in S_n(\bZ):(\zs A\in\cF)x\in S_n(A)\}$. 

The monad of $\cF$ is $\mu(\cF)=\bigcup_{n<\omega}\mu_n(\cF)=\{x\in\tac \bZ:(\zs A\in\cF)x\in\tac A\}$.

For $x\in\tac \bZ$, $v(x)$ is the unique $\cF\in\beta\bZ$ such that $x\in\mu(\cF)$.
\end{de}

Note that this definition of $v(x)$ agrees with the previous one (for $x\in\zve \bZ$).

\begin{pp}\label{monadi}
(\cite{L1}, Proposition 2.5.11) For every $x\in\tac \bZ$ and every $n\in\omega$, $v(S_n(x))=v(x)$.
\end{pp}

Let us recall the tensor (or Fubini) product of ultrafilters: for $\cF,\cG\in\beta \bZ$, $\cF\otimes\cG$ is the ultrafilter on $\bZ\times \bZ$ defined by
$$S\in\cF\otimes\cG\dl\{x\in \bZ:\{y\in \bZ:(x,y)\in S\}\in\cG\}\in\cF.$$
The definitions of monads of ultrafilters of the form $\cF\otimes\cG$ and the corresponding function $v$ are analogous as above. For ultrafilters $\cF$ and $\cG$ and nonstandard numbers $x\in\mu(\cF)$ and $y\in\mu(\cG)$, $(x,y)$ is a {\it tensor pair} if $(x,y)\in\mu(\cF\otimes\cG)$. 

\begin{lm}\label{minustensor}
If $(x,y)\in\zve\bZ\times\zve\bZ$ is a tensor pair, then so are $(x,-y)$ and $(-x,y)$.
\end{lm}

\dokaz Let $\cF=v(x)$ and $\cG=v(y)$; then $v(-y)=-\cG$ and $v((x,y))=\cF\otimes\cG$. We need to prove that $v((x,-y))=\cF\otimes(-\cG)$. But whenever $(x,-y)\in\zve S$ for some $S\subseteq\bZ\times\bZ$, we have $(x,y)\in\zve S'$, where $S':=\{(m,-n):(m,n)\in S\}$. By the assumptions $S'\in\cF\otimes\cG$, so $\{x\in\bZ:\{y\in\bZ:(x,y)\in S\}\in(-\cG)\}=\{x\in\bZ:-\{y\in\bZ:(x,y)\in S\}\in\cG\}=\{x\in\bZ:\{y\in\bZ:(x,y)\in S'\}\in\cG\}\in\cF$, so $S\in\cF\otimes(-\cG)$.

The proof for $(-x,y)$ is analogous.\kraj

By \cite{DN}, Proposition 11.7.2, for any tensor pair $(x,y)$ we have $x+y\in\mu(\cF+\cG)$ and $x\cdot y\in\mu(\cF\cdot\cG)$. An important feature of $\omega$-hyperextensions is that they provide a canonical way to obtain tensor pairs.

\begin{pp}\label{zbirpro}
(\cite{L1}, Theorem 2.5.27) If $x\in\mu(\cF)$ and $y\in\mu(\cG)$, then the pair $(x,\zve y)$ is a tensor pair. Hence, $x+\zve y\in\mu(\cF+\cG)$ and $x\cdot\zve y\in\mu(\cF\cdot\cG)$.
\end{pp}

\section{Congruence modulo ultrafilter}\label{flawed}

A natural way to define the congruence relation modulo an ultrafilter would be to imitate again the construction of an extension $\widetilde{\rho}$, as described in Section \ref{kongm}. 

\begin{de}\label{defcong}
For $\cM\in\beta \bN$ and $\cF,\cG\in\beta \bZ$, $\cF\equiv_\cM\cG$ if and only if for every $A\in\cM$ the set $\{(x,y)\in \bZ\times \bZ:(\po m\in A)x\equiv_my\}$ belongs to the ultrafilter $\cF\otimes\cG$.
\end{de}

This definition has a nice equivalent formulation via divisibility of ultrafilters.

\begin{lm}\label{equivdelj}
For $\cM\in\beta \bN$ and $\cF,\cG\in\beta \bZ$, $\cF\equiv_\cM\cG$ if and only if $\cM\widemid\cF-\cG$.
\end{lm}

\dokaz \begin{eqnarray*}
\cF\equiv_\cM\cG &\dl & (\zs A\in\cM)\{x\in \bZ:\{y\in \bZ:(\po m\in A)x\equiv_my\}\in\cG\}\in\cF\\
  &\dl & (\zs A\in\cM)\{x\in \bZ:\{y\in \bZ:x-y\in A\gstr\}\in\cG\}\in\cF\\
  &\dl & (\zs A\in\cM\cap\cU)\{x\in \bZ:\{y\in \bZ:x-y\in A\}\in\cG\}\in\cF\\
  &\dl & (\zs A\in\cM\cap\cU)\{x\in \bZ:x-A\in\cG\}\in\cF\\
  &\dl & (\zs A\in\cM\cap\cU)A\in\cF-\cG,
 \end{eqnarray*}
which is equivalent to $\cM\widemid\cF-\cG$.\kraj

The following lemma justifies our using the same notation as for the relation from Section \ref{kongm}.

\begin{lm}\label{restr1}
If $m\in \bN$ and $\cF,\cG\in\beta \bZ$, $\cF\equiv_m\cG$ as defined in Section \ref{kongm} is equivalent to $\cF\equiv_m\cG$ from Definition \ref{defcong}.
\end{lm}

\dokaz Since $h_m$ is a homomorphism, $\widetilde{h_m}(\cF-\cG)=\widetilde{h_m}(\cF)-\widetilde{h_m}(\cG)$. It follows that $m\mid\cF-\cG$ if and only if $\widetilde{h_m}(\cF-\cG)=0$, if and only if $\widetilde{h_m}(\cF)-\widetilde{h_m}(\cG)$.\kraj


$\equiv_\cM$ also has a nonstandard characterization. First we recall Puritz's result that $(x,y)\in\zve\bN\times\zve\bN$ is a tensor pair if and only if $x<\zve f(y)$ for every $f:\bN\str \bN$ such that $\zve f(y)\in\zve \bN\setminus \bN$ (\cite{P2}, Theorem 3.4). Taking into account Lemma \ref{minustensor}, we get the following version of this result.

\begin{pp}
$(x,y)\in\zve\bZ\times\zve\bZ$ is a tensor pair if and only if $|x|<|\zve f(y)|$ for every $f:\bZ\str\bZ$ such that $\zve f(y)\in\zve\bZ\setminus\bZ$.
\end{pp}

If we denote $\cG=v(y)$, the condition $\zve f(y)\notin\bZ$ is equivalent to $f\rest B$ not being constant for any $B\in\cG$. Let us call $f:\bZ\str \bZ$ non-$\cG$-constant in that case.

Note that we are still working in any ${\goth c}^+$-enlargement (we do not need an $\omega$-hyperextension), so $\mu(\cF)$ here actually means $\mu_1(\cF)$.

\begin{te}
Let $\cM\in\beta\bN$ and $\cF,\cG\in\beta\bZ$. The following conditions are equivalent:

(i) $\cF\equiv_\cM\cG$;

(ii) in some ${\goth c}^+$-enlargement holds
\begin{equation}\label{eqpair}
(\zs m\in\mu(\cM))(\po x\in\mu(\cF))(\po y\in\mu(\cG))((x,y)\mbox{ is a tensor pair }\land m\mid x-y)
\end{equation}

(iii) in every ${\goth c}^+$-enlargement holds (\ref{eqpair}).
\end{te}

\dokaz (ii)$\Str$(i) Let (\ref{eqpair}) hold in some ${\goth c}^+$-enlargement. If $y\in\mu(\cG)$ then $-y\in\mu(-\cG)$. Since for a tensor pair $(x,y)$ we have, by Lemma \ref{minustensor}, $x-y=x+(-y)\in\mu(\cF-\cG)$, the ``if" part follows directly from Proposition \ref{ekviv}.\\

(i)$\Str$(iii) Assume $\cM\widemid\cF-\cG$; we work in arbitrary ${\goth c}^+$-enlargement. We define, for $A,B\subseteq\bZ$, $M\subseteq \bN$ and $f:\bZ\str\bZ$:
\begin{eqnarray*}
& & E_{A,B,M}=\{(m,a,b)\in\bN\times\bZ\times\bZ:a\in A\land b\in B\land m\in M\land m\mid a-b\}\\
& & F_f=\{(m,a,b)\in\bN\times\bZ\times\bZ:|a|<|f(b)|\}.
\end{eqnarray*}
We prove that the family $\{E_{A,B,M}:A\in\cF,B\in\cG,M\in\cM\}\cup\{F_f:f:\bZ\str\bZ\mbox{ is non-}\cG\mbox{-constant}\}$ has the finite intersection property. $\{E_{A,B,M}:A\in\cF,B\in\cG,M\in\cM\}$ is closed for finite intersections. So let $A\in\cF$, $B\in\cG$, $M\in\cM$ and let $f_1,f_2,\dots,f_k:\bZ\str\bZ$ be non-$\cG$-constant. Since $M\gstr\in\cM\cap\cU$, $\cM\widemid\cF-\cG$ implies $M\gstr\in\cF-\cG$. Hence $\{n\in\bZ:n-M\gstr\in\cG\}\in\cF$. Let $a\in A\cap\{n\in\bZ:n-M\gstr\in\cG\}$. This means that $B_1:=B\cap(a-M\gstr)\in\cG$. Hence there is $b\in B_1$ such that $|f_i(b)|>|a|$ for all $i\leq k$ (otherwise $\{b\in B_1:f_i(b)=j\}\in\cG$ for some $i\leq k$ and some $-a\leq j\leq a$, a contradiction with the assumption that $f_i$ is non-$\cG$-constant). Since $b\in a-M\gstr$, there is $m\in M$ such that $m\mid a-b$, so $(m,a,b)\in E_{A,B,M}\cap F_{f_1}\cap F_{f_2}\cap\dots\cap F_{f_k}$.

Now, since we are working with a ${\goth c}^+$-enlargement, there is
$$(m,x,y)\in\bigcap_{A\in\cF,B\in\cG,M\in\cM}\zve E_{A,B,M}\;\;\cap\bigcap_{f\mbox{ non-}\cG\mbox{-constant}}\zve F_f.$$
This means that $m\in\mu(\cM)$, $x\in\mu(\cF)$, $y\in\mu(\cG)$ and $m\mid x-y$. Also, for every non-$\cG$-constant $f:\bZ\str\bZ$, $|\zve f(y)|>|x|$, so $(x,y)$ is a tensor pair.\kraj

Unfortunately, we do not even know whether $\equiv_\cM$ is an equivalence relation on $\beta \bZ$, which makes it unconvenient to work with. Therefore in the next section we introduce a stronger relation with much nicer properties.

\section{Strong congruence}\label{congult}


To better explain the forthcoming definition of congruence, we begin with a few simple lemmas. Recall that $MAX$ is the class of ultrafilters $\widemid$-divisible by all others.

\begin{lm}\label{razlika1}
Let $x,y\in\tac \bZ$ and $v(x)=v(y)$. Then $m\mid x-y$ for all $m\in \bN$ and $x-y\in\mu(MAX)$.
\end{lm}

\dokaz For each $m\in \bN$, let $h_m$ be the function defined in Section \ref{kongm}. Then $\tac h_m(x)\in \bZ_m$ for all $x\in\tac \bZ$. By Proposition \ref{slaganjev}, $v(\tac h_m(x))=\widetilde{h_m}(v(x))=\widetilde{h_m}(v(y))=v(\tac h_m(y))$, so $x$ and $y$ have the same residue modulo $m$.

Ultrafilters from $MAX$ are those divisible by all $m\in\bN$. Hence $\mu(MAX)$ consists exactly of nonstandard numbers divisible by all $m\in\bN$, so the second statement follows directly from the first.\kraj

By Theorem \ref{geomset}, the assumption of Lemma \ref{razlika1} can not be relaxed to $v(x)=_\sim v(y)$: there are $=_\sim$-equivalent ultrafilters giving different residues modulo some $m\in \bN$.

\begin{lm}\label{razlika2}
Let $x,y\in\tac \bZ$, $v(x)=v(y)$ and $m\in S_k(\bN)$. Then $m\mid S_k(x)-S_k(y)$.
\end{lm}

\dokaz By Lemma \ref{razlika1}, $(\zs m\in \bN)m\mid x-y$. By Transfer, $(\zs m\in S_k(\bN))m\mid S_k(x)-S_k(y)$.\kraj

Thus, for every $m\in S_k(\bN)$, all the numbers from $\mu(\cF)\cap S_k[\tac \bZ]$ have the same residue modulo $m$. We will use this to establish a strengthening of congruence modulo $\cM\in\beta \bN$.

\begin{de}\label{defkong}
Ultrafilters $\cF,\cG\in\beta \bZ$ are strongly congruent modulo $\cM\in\beta \bN$ if, in every $\omega$-hyperenlargement,
\begin{equation}\label{eq3}
(\zs m\in\mu_1(\cM))(\po x\in\mu(\cF))(\po y\in\mu(\cG))m\mid\zve x-\zve y.
\end{equation}
We write $\cF\equiv_\cM^s\cG$.
\end{de}

We easily get the following equivalent condition.

\begin{lm}\label{strongcong}
$\cF\equiv_\cM^s\cG$ implies that in every $\omega$-hyperenlargement
$$(\zs m\in\mu_1(\cM))(\zs x\in\mu(\cF))(\zs y\in\mu(\cG))m\mid\zve x-\zve y.$$
\end{lm}

\dokaz Let $x_0\in\mu(\cF)$ and $y_0\in\mu(\cG)$ be such that $m\mid\zve x_0-\zve y_0$, and let $x\in\mu(\cF)$ and $y\in\mu(\cG)$ be arbitrary. By Lemma \ref{razlika2}, $m\mid\zve x-\zve x_0$ and $m\mid\zve y-\zve y_0$, so $m\mid\zve x-\zve y$ as well.\kraj

To avoid constant repetition, in each of the proofs in the rest of the paper it will be understood that we are working in an $\omega$-hyperenlargement (a single structure extension which is a ${\goth c}^+$-enlargement).

It will follow from Lemmas \ref{strongovi}, \ref{deljivosti} and \ref{equivdelj} that $\cF\equiv_\cM^s\cG$ implies $\cF\equiv_\cM\cG$. For now we prove that $\equiv_m^s$ for $m\in\bN$ also coincides with the congruence relation modulo integer (from Section \ref{kongm}).

\begin{lm}\label{restr}
If $m\in\bN$ and $\cF,\cG\in\beta \bZ$, $\cF\equiv_m^s\cG$ holds if and only if $\cF\equiv_m\cG$.
\end{lm}

\dokaz The only element of $\mu_1(m)$ is $m$ itself. Let $x\in\mu(\cF)$ and $y\in\mu(\cG)$ be such that $m\mid\zve x-\zve y$; then $\zve x$ and $\zve y$ have the same residue modulo $m$: $\tac h_m(\zve x)=\tac h_m(\zve y)$. Then, by Propositions \ref{slaganjev} and \ref{monadi}, $\widetilde{h_m}(\cF)=v(\tac h_m(\zve x))=v(\tac h_m(\zve y))=\widetilde{h_m}(\cG)$, so $\cF\equiv_m\cG$. The other implication is proved similarly, using Lemma \ref{strongcong}.\kraj

\begin{lm}
$\equiv_\cM^s$ is an equivalence relation on the set $\beta \bZ$.
\end{lm}

\dokaz Reflexivity and symmetry are obvious from the definition. So let $\cF\equiv_\cM^s\cG$ and $\cG\equiv_\cM^s\cH$. By Lemma \ref{strongcong}, for any $m\in\mu_1(\cM)$, $x\in\mu(\cF)$, $y\in\mu(\cG)$ and $z\in\mu(\cH)$ holds $m\mid\zve x-\zve y$ and $m\mid\zve y-\zve z$. Then $m\mid\zve x-\zve z$, so $\cF\equiv_\cM^s\cH$.\kraj

\begin{te}\label{compat}
Let $\cM\in\beta \bN$. $\equiv_\cM^s$ is compatible with operations $+$ and $\cdot$ in $\beta \bZ$:

(a) $\cF_1\equiv_\cM^s\cF_2$ and $\cG_1\equiv_\cM^s\cG_2$ imply $\cF_1+\cG_1\equiv_\cM^s\cF_2+\cG_2$;

(b) $\cF_1\equiv_\cM^s\cF_2$ and $\cG_1\equiv_\cM^s\cG_2$ imply $\cF_1\cdot\cG_1\equiv_\cM^s\cF_2\cdot\cG_2$.
\end{te}

\dokaz Let $m\in\mu_1(\cM)$, $x_1\in\mu_1(\cF_1)$, $x_2\in\mu_1(\cF_2)$, $y_1\in\mu_1(\cG_1)$ and $y_2\in\mu_1(\cG_2)$. It follows from Proposition \ref{monadi} that $\zve y_1\in\mu(\cG_1)$ and $\zve y_2\in\mu(\cG_2)$. By the assumptions we have $m\mid\zve x_1-\zve x_2$ and $m\mid\zve{\zve y_1}-\zve{\zve y_2}$.

(a) By Proposition \ref{zbirpro} $x_1+\zve y_1\in\mu(\cF_1+\cG_1)$ and $x_2+\zve y_2\in\mu(\cF_2+\cG_2)$. From the above conclusions follows $m\mid(\zve x_1+\zve{\zve y_1})-(\zve x_2+\zve{\zve y_2})$, i.e.\ $m\mid\zve(x_1+\zve y_1)-\zve(x_2+\zve y_2)$. Since we started with arbitrary $m\in\mu_1(\cM)$, this means that $\cF_1+\cG_1\equiv_\cM^s\cF_2+\cG_2$.

(b) By Proposition \ref{zbirpro} $x_1\cdot\zve y_1\in\mu(\cF_1\cdot\cG_1)$ and $x_2\cdot\zve y_2\in\mu(\cF_2\cdot\cG_2)$. We have $m\mid(\zve x_1-\zve x_2)\zve{\zve y_1}$ and $m\mid\zve x_2(\zve{\zve y_1}-\zve{\zve y_2})$. Hence $m\mid\zve x_1\zve{\zve y_1}-\zve x_2\zve{\zve y_2}$, i.e.\ $m\mid\zve(x_1\zve y_1)-\zve(x_2\zve y_2)$, so $\cF_1\cdot\cG_1\equiv_\cM^s\cF_2\cdot\cG_2$.\kraj

%




The following simple result is a version of a well-known fact (\cite{P1}, Corollary 8.3).

\begin{lm}\label{FminusF}
(a) Every $\cF\in MAX$ is strongly congruent to zero modulo any ultrafilter;

(b) for every $\cF\in\beta \bZ\setminus \bZ$, $\cF-\cF\in MAX$.
\end{lm}

\dokaz (a) For any $\cF\in MAX$ and any $x\in\mu(\cF)$, $(\zs m\in \bN)m\mid x$ implies by Transfer $(\zs m\in\zve \bN)m\mid\zve x$, which gives us $\cF\equiv_\cM 0$ for any $\cM$.\\

(b) We will show that $A\in\cF-\cF$ for all $A\in\cU_Z$. Let $m\in A$ be arbitrary. Then there is $r\in \bZ_m$ such that $m\bZ+r\in\cF$, so since $m\bZ\subseteq -A$, it follows that $n-A\in\cF$ for all $n\in m\bZ+r$. Thus $m\bZ+r\subseteq\{n\in \bZ:n-A\in\cF\}$, so $\{n\in \bZ:n-A\in\cF\}\in\cF$, which means that $A\in \cF-\cF$.\kraj

Let us also note, regarding the lemma above, that $\cF=_\sim\cG$ is not enough to conclude that $\cF-\cG\in MAX$. By Theorem \ref{geomset} there are $\cF,\cG\in\beta \bN$ and $m\in \bN$ such that $\cF=_\sim\cG$ but $\cF\not\equiv_m\cG$, say $\cF\equiv_mr_1$ and $\cG\equiv_mr_2$ for some $r_1<m$ and $r_2<m$. From Proposition \ref{hom} we get $\cF-\cG\equiv_m r_1-r_2\neq 0$, so $m\nmid\cF-\cG$.

\begin{de}
A family $\{\cF_i:i\in I\}$ of ultrafilters is a complete residue system modulo $\cM\in\beta \bN$ if it contains exactly one element of every equivalence class of strong congruence modulo $\cM$.
\end{de}

As an application of the above results, we have an ultrafilter version of a well-known theorem on complete residue systems in $\bZ$.

\begin{te}\label{crs}
If $\{\cF_i:i\in I\}$ is a complete residue system modulo $\cM\in\beta \bN$ then, for every $\cG\in\beta \bN$, $\{\cF_i+\cG:i\in I\}$ and $\{\cG+\cF_i:i\in I\}$ are complete residue systems modulo $\cM$.
\end{te}

\dokaz We need to show that in $\cR=\{\cF_i+\cG:i\in I\}$ no two ultrafilters are congruent modulo $\cM$, and that each congruence class has a representative in $\cR$.

First assume $\cF_i+\cG\equiv_\cM^s\cF_j+\cG$ for some $i,j\in I$, $i\neq j$. By Theorem \ref{compat} $\cF_i+\cG-\cG\equiv_\cM^s\cF_j+\cG-\cG$. By Lemma \ref{FminusF} $\cF_i=\cF_i+0\equiv_\cM^s\cF_i+\cG-\cG\equiv_\cM^s\cF_j+\cG-\cG\equiv_\cM^s\cF_j$, a contradiction.

Now let $\cH\in\beta \bN$ be arbitrary. There is $i\in I$ such that $\cF_i\equiv_\cM^s\cH-\cG$. Using Theorem \ref{compat} and Lemma \ref{FminusF} again we get $\cF_i+\cG\equiv_\cM^s\cH-\cG+\cG\equiv_\cM^s\cH$.

The proof that $\{\cG+\cF_i:i\in I\}$ is a complete residue system modulo $\cM\in\beta \bN$ is analogous.\kraj

\section{Strong divisibility}

It is natural to ask: which ultrafilters are strongly congruent to zero modulo some $\cM\in\beta \bN$? Are those exactly the ultrafilters divisible by $\cM$? For example, we saw in Lemma \ref{FminusF} that $\widemid$-maximal ultrafilters are always strongly congruent to zero. In general, the above question leads us to the following definition.

\begin{de}\label{defstrong}
Let $\cM\in\beta \bN$ and $\cF\in\beta \bZ$. $\cF$ is strongly divisible by $\cM$ if, in every $\omega$-hyperenlargement,
$$(\zs m\in\mu_1(\cM))(\po x\in\mu(\cF))m\mid\zve x.$$
We write $\cM\mid^s\cF$.
\end{de}

In the same way as Lemma \ref{strongcong}, we get a seemingly stronger condition.

\begin{lm}\label{strongstrong}
$\cM\mid^s\cF$ implies that in every $\omega$-hyperenlargement
$$(\zs m\in\mu_1(\cM))(\zs x\in\mu(\cF))m\mid\zve x.$$
\end{lm}

Proposition \ref{ekviv} easily implies the following.

\begin{lm}\label{deljivosti}
For all $\cM\in\beta \bN$ and $\cF\in\beta \bZ$, $\cM\mid^s\cF$ implies $\cM\widemid\cF$.
\end{lm}


It is tempting to try to prove the reverse implication; unfortunately this is not true, as we will now see.

\begin{lm}\label{freeprost}
No $\bN$-free ultrafilter has any $\mid^s$-divisors.
\end{lm}

\dokaz Assume the opposite, that an $\bN$-free ultrafilter $\cF$ is $\mid^s$-divisible by some $\cG$. Then $\cG$ is also $\bN$-free, so for any $x\in\mu(\cF)$ holds $(\zs m\in \bN)m\nmid x$. By Transfer $(\zs m\in\zve \bN)m\nmid\zve x$, a contradiction with $\cG\nmid^s\cF$.\kraj

Thus, this notion of divisibility is too strong to be our main divisibility relation, but it has some properties that are in good accordance with the strong congruence relation and operations on $\beta \bN$.

However, Lemma \ref{freeprost} also says that $\mid^s$ is not reflexive: $\bN$-free ultrafilters are not divisible by themselves. It is, however, transitive: let $\cF\mid^s\cG$ and $\cG\mid^s\cH$. Let $x\in\mu_1(\cF)$, $y\in\mu_1(\cG)$ and $z\in\mu_1(\cH)$ be arbitrary. Then $x\mid\zve y$ and $y\mid\zve z$. Hence $\zve y\mid\zve{\zve z}$, so $x\mid\zve{\zve z}$, which suffices for $\cF\mid^s\cH$.

\begin{lm}\label{strongovi}
$\cF\equiv_\cM^s\cG$ if and only if $\cM\mid^s\cF-\cG$.
\end{lm}

\dokaz ($\Rightarrow$) Let $m\in\mu_1(\cM)$ be arbitrary and let $x\in\mu_1(\cF)$ and $y\in\mu_1(\cG)$ be such that $m\mid\zve x-\zve y$. By Proposition \ref{monadi}, $v(y)=v(\zve y)$ so, by Lemma \ref{razlika2}, $m\mid\zve y-\zve{\zve y}$. It follows that $m\mid\zve x-\zve{\zve y}$, i.e.\ $m\mid\zve(x-\zve y)$. On the other hand, since $-y\in\mu(-\cG)$, by Lemma \ref{minustensor} and Proposition \ref{zbirpro}, $x-\zve y=x+\zve(-y)\in\mu(\cF-\cG)$, so $\cM\mid^s\cF-\cG$.

($\Leftarrow$) Let $m\in\mu_1(\cM)$, $x\in\mu_1(\cF)$ and $y\in\mu_1(\cG)$ be arbitrary. Then $x-\zve y\in\mu(\cF-\cG)$ so, by Lemma \ref{strongstrong}, $m\mid\zve(x-\zve y)$. By Lemma \ref{razlika2} again we have $m\mid\zve y-\zve{\zve y}$, so $m\mid\zve x-\zve y$, meaning that $\cF\equiv_\cM^s\cG$.\kraj

\begin{te}
Let $\cM\in\beta \bN$ and $\cF,\cG\in\beta \bZ$.


(a) $\cM\mid^s\cF$ and $\cM\mid^s\cG$ imply $\cM\mid^s\cF+\cG$;

(b) $\cM\mid^s\cF$ implies $\cM\mid^s\cF\cdot\cG$;

(c) $\cM\mid^s\cG$ implies $\cM\mid^s\cF\cdot\cG$.
\end{te}

\dokaz Let $m\in\mu_1(\cM)$, $x\in\mu_1(\cF)$ and $y\in\mu_1(\cG)$.


(a) By assumptions $m\mid\zve x$ and $m\mid\zve{\zve y}$. Hence $m\mid\zve(x+\zve y)$, and therefore $\cM\mid^s\cF+\cG$.

(b) Now we have $m\mid\zve x$, which suffices for $m\mid\zve x\zve{\zve y}$ i.e.\ $m\mid\zve(x\zve y)$, so $\cM\mid^s\cF\cdot\cG$.

(c) By Lemma \ref{strongstrong} $\cM\mid^s\cG$ implies $m\mid\zve{\zve y}$, so again $m\mid\zve x\zve{\zve y}$ and $\cM\mid^s\cF\cdot\cG$.\kraj

Let us remind ourselves of the definitions of other three divisibility relations from \cite{So1}:
\begin{eqnarray*}
\cG\mid_L\cF & \mbox{iff} & (\po\cH\in\beta \bN)\cF=\cH\cdot\cG\\
\cG\mid_R\cF & \mbox{iff} & (\po\cH\in\beta \bN)\cF=\cG\cdot\cH\\
\cG\mid_M\cF & \mbox{iff} & (\po\cH_1,\cH_2\in\beta \bN)\cF=\cH_1\cdot\cG\cdot\cH_2.
\end{eqnarray*}

What is the place of $\mid^s$ (restricted to $\beta\bN\times\beta\bN$) among these relations? Like all the others, its restriction to $\bN\times\bN$ is just the usual divisibility relation (Lemma \ref{restr}). We already saw that $\mid^s\subset\widemid$. We will show that this is the only inclusion that can be established:

\begin{center}
\includegraphics[scale=0.15]{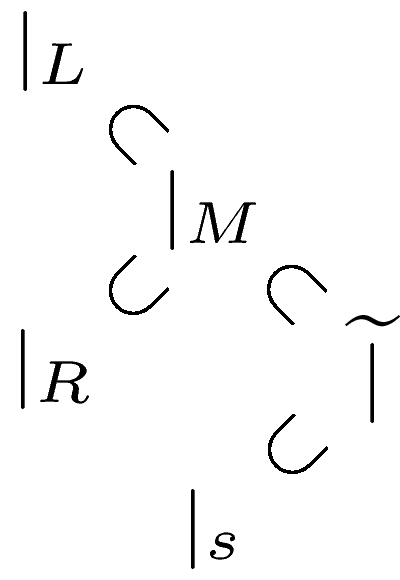}
\end{center}

First, why $\mid_L\not\subseteq\mid^s$? Let $\cP,\cQ\in\beta \bN\setminus \bN$ be $\widemid$-prime and let $\cF=\cP\cdot\cQ$. Then $\cQ\mid_L\cF$ but, by Lemma \ref{freeprost}, $\cQ\nmid_s\cF$. Analogously we conclude that $\mid_R\not\subseteq\mid^s$.

That $\mid^s\subseteq\mid_M$ does not hold either can be seen by considering maximal classes of these two orders. By \cite{So2}, Theorem 4.1, the $\mid_M$-maximal ultrafilters are exactly those in the smallest ideal $K(\beta \bN,\cdot)$. On the other hand, the class of $\mid^s$-maximal ultrafilters is exactly $MAX$ by Lemmas \ref{FminusF} and \ref{deljivosti}. But $MAX$ is a proper superset of $K(\beta \bN,\cdot)$; we postpone the detailed examination of this and other aspects of maximal ultrafilters until a projected sequel to this paper.

\section{Final remarks and questions}





Even after finding, in Section \ref{flawed}, several equivalent conditions for $\equiv_\cM$, we were not able to answer the following.

\begin{qu}
Is $\equiv_\cM$ an equivalence relation?
\end{qu}

Not being able to prove that it is presents a big drawback for using this relation, which seems to be the most natural extension of the congruence relation to $\beta \bN$.

Some more properties of our relations could be proved if we worked with ${\goth c}^+$-saturated nonstandard extensions. This is a stronger condition than being a ${\goth c}^+$-enlargement: $(V(Y),*)$ is $\kappa$-{\it saturated} if every family $F$ of internal sets in $V(Y)$ with the finite intersection property such that $|F|<\kappa$ has nonempty intersection. To Proposition \ref{ekviv} one can add two more equivalent conditions (see \cite{So5}, Theorem 3.4):\\

(iv) in every ${\goth c}^+$-saturated extension $V(Y)$, for every $x\in\mu(\cF)$ there is $y\in\mu(\cG)$ such that $x\zvez\mid y$;

(v) in every ${\goth c}^+$-saturated extension $V(Y)$, for every $y\in\mu(\cG)$ there is $x\in\mu(\cF)$ such that $x\zvez\mid y$.\\


However, Proposition \ref{prenosnatac} does not hold for ${\goth c}^+$-saturation in place of ${\goth c}^+$-enlargement: see \cite{L1}, page 74. So to use the equivalents (iv) and (v) we would have to answer the following question.


\begin{qu}
Is it possible to construct a ${\goth c}^+$-saturated $\omega$-hyperextension of $\bZ$?
\end{qu}

\section{Declarations}

The author acknowledges financial support of the Science Fund of the Republic of Serbia (call PROMIS, project CLOUDS, grant no.\ 6062228) and Ministry of Education, Science and Technological Development of the Republic of Serbia (grant no.\ 451-03-68/2020-14/200125).





\end{document}